\documentclass[10pt,twoside]{article}
\usepackage{amsmath,epsf,amssymb,cite}
\def\firstpage{1}\def\lastpage{1000}

\makeatletter
\def\magnification{\afterassignment\m@g\count@}
\def\m@g{\mag=\count@\hsize6.5truein\vsize8.9truein\dimen\footins8truein}
\makeatother

\DeclareFontFamily{U}{rsf}{}
\DeclareFontShape{U}{rsf}{m}{n}{
  <5> <6> rsfs5 <7> <8> <9> rsfs7 <10-> rsfs10}{}
\DeclareMathAlphabet\Scr{U}{rsf}{m}{n}

\font\caps=cmcsc10                    
\font\Caps=cccsc10 scaled \magstep1   

\makeatletter

\renewcommand\section{\@startsection {section}{1}{\z@}%
                                   {-3.5ex \@plus -1ex \@minus -.2ex}%
                                   {2.3ex \@plus.2ex}%
                                   {\normalfont\Large\normalsize\caps}}

\@specialpagefalse\headheight=8.5pt
\def\pplogo{\vbox{\kern-\headheight\kern -29pt
\halign{##&##\hfil\cr&{
\ppnumber}\cr\rule{0pt}{2.5ex}&\ppdate\cr}}}
\def\ppnumber{\vbox{\baselineskip14pt\hbox{DUKE-CGTP-98-06}
\hbox{math.AG/9809004}}}
\def\ppdate{June 1998} \date{}
\def\DocMath{\pplogo}
\renewcommand{\@oddfoot}{}
\renewcommand{\@evenfoot}{}%
\renewcommand{\@evenhead}{%
    \ifnum\thepage>\lastpage\rlap{\thepage}\hfill%
    \else\rlap{\thepage}\slshape\leftmark\hfill\caps\SAuthor\hfill\fi}%
\renewcommand{\@oddhead}{%
    \ifnum\thepage=\firstpage{\DocMath\hfill\llap{\thepage}}%
    \else{\slshape\rightmark}\hfill\caps\STitle\hfill\llap{\thepage}\fi}%
\makeatother

\def\TSkip{\bigskip}
\newbox\TheTitle{\obeylines\gdef\GetTitle #1
\ShortTitle  #2
\SubTitle    #3
\Author      #4
\ShortAuthor #5
\EndTitle
{\setbox\TheTitle=\vbox{\baselineskip=20pt\let\par=\cr\obeylines%
\halign{\centerline{\Caps##}\cr\noalign{\medskip}\cr#1\cr}}%
	\copy\TheTitle\TSkip\TSkip%
\def\next{#2}\ifx\next\empty\gdef\STitle{#1}\else\gdef\STitle{#2}\fi%
\def\next{#3}\ifx\next\empty%
    \else\setbox\TheTitle=\vbox{\baselineskip=20pt\let\par=\cr\obeylines%
    \halign{\centerline{\caps##} #3\cr}}\copy\TheTitle\TSkip\TSkip\fi%
\centerline{\caps #4}\TSkip\TSkip%
\def\next{#5}\ifx\next\empty\gdef\SAuthor{#4}\else\gdef\SAuthor{#5}\fi%
\catcode'015=5}}\def\Title{\obeylines\GetTitle}
\def\Abstract{\begingroup\narrower
    \parskip=\medskipamount\parindent=0pt{\caps Abstract. }}
\def\EndAbstract{\par\endgroup\TSkip}

\long\def\MSC#1\EndMSC{\def\arg{#1}\ifx\arg\empty\relax\else
     {\par\narrower\noindent%
     1991 Mathematics Subject Classification: #1\par}\fi}

\long\def\KEY#1\EndKEY{\def\arg{#1}\ifx\arg\empty\relax\else
	{\par\narrower\noindent Keywords and Phrases: #1\par}\fi\TSkip}

\newbox\TheAdd\def\Addresses{\vfill\copy\TheAdd\vfill
    \ifodd\number\lastpage\vfill\eject\phantom{.}\vfill\eject\fi}
{\obeylines\gdef\GetAddress #1
\Address #2
\Address #3
\Address #4
\EndAddress
{\def\xs{6truecm}
\setbox0=\vtop{{\obeylines\hsize=\xs#1\par}}\def\next{#2}
\ifx\next\empty 
     \setbox\TheAdd=\hbox to\hsize{\hfill\copy0\hfill}
\else\setbox1=\vtop{{\obeylines\hsize=\xs#2\par}}\def\next{#3}
\ifx\next\empty 
     \setbox\TheAdd=\hbox to\hsize{\hfill\copy0\hfill\copy1\hfill}
\else\setbox2=\vtop{{\obeylines\hsize=\xs#3\par}}\def\next{#4}
\ifx\next\empty\ 
     \setbox\TheAdd=\vtop{\hbox to\hsize{\hfill\copy0\hfill\copy1\hfill}
                \vskip20pt\hbox to\hsize{\hfill\copy2\hfill}}
\else\setbox3=\vtop{{\obeylines\hsize=\xs#4\par}}
     \setbox\TheAdd=\vtop{\hbox to\hsize{\hfill\copy0\hfill\copy1\hfill}
	        \vskip20pt\hbox to\hsize{\hfill\copy2\hfill\copy3\hfill}}
\fi\fi\fi\catcode'015=5}}\gdef\Address{\obeylines\GetAddress}

\begin{document}

\def\O{{\mathcal O}}
\def\C{{\mathbb C}}
\def\P{{\mathbb P}}
\def\Q{{\mathbb Q}}
\def\R{{\mathbb R}}
\def\Z{{\mathbb Z}}
\def\H{{\mathbb H}}
\def\M{\Scr{M}}
\def\mQ{\Scr{Q}}
\def\Hom{\operatorname{Hom}}
\def\Ext{\operatorname{Ext}}
\def\Tors{\operatorname{Tors}}
\def\Free{\operatorname{Free}}
\def\Ker{\operatorname{Ker}}
\def\Spec{\operatorname{Spec}}
\def\Area{\operatorname{Area}}
\def\Vol{\operatorname{Vol}}
\def\ad{\operatorname{ad}}
\def\ch{\operatorname{ch}}
\def\td{\operatorname{td}}
\def\tr{\operatorname{tr}}
\def\Pic{\operatorname{Pic}}
\def\Jac{\operatorname{Jac}}
\def\disc{\operatorname{disc}}
\def\Img{\operatorname{Im}}
\def\Rea{\operatorname{Re}}
\def\Gr{\operatorname{Gr}}
\def\SO{\operatorname{SO}}
\def\Sl{\operatorname{SL}}
\def\GO{\operatorname{O{}}}
\def\SU{\operatorname{SU}}
\def\GU{\operatorname{U{}}}
\def\Sp{\operatorname{Sp}}
\def\rank{\operatorname{rank}}
\def\Spin{\operatorname{Spin}}
\def\so{\operatorname{\mathfrak{so}}}
\def\su{\operatorname{\mathfrak{su}}}
\def\gu{\operatorname{\mathfrak{u}}}
\def\sp{\operatorname{\mathfrak{sp}}}
\def\labto#1{\mathrel{\mathop\to^{#1}}}
\def\Dr{\mathbf{D}}
\def\sm{$\sigma$-model}
\def\nlsm{non-linear \sm}
\def\MW{Mordell--Weil}
\def\CY{Calabi--Yau}
\def\ff#1#2{{\textstyle\frac{#1}{#2}}}
\def\RoR{$R\leftrightarrow1/R$}
\def\spnh{\Spin(32)/\Z_2}
\def\mbf#1{\mbox{\mathversion{bold}$#1$}}
\def\HS#1{{\mathbb{F}}_{#1}}

\newtheorem{theorem}{Theorem}
\newtheorem{prop}{Proposition}
\newtheorem{defn}{Definition}

\Title
String Theory and Duality
\ShortTitle
\SubTitle
\Author
Paul S. Aspinwall
\ShortAuthor
\EndTitle
\Abstract
String Duality is the statement that one kind of string theory
compactified on one space is equivalent in some sense to another
string theory compactified on a second space. This draws a connection
between two quite different spaces. Mirror symmetry is an example of
this. Here we discuss mirror symmetry and another ``heterotic/type
II'' duality  which relates vector bundles on
a K3 surface to a \CY\ threefold.
\EndAbstract
\MSC
81T30, 14J32, 14J28, 14J60
\EndMSC
\KEY
\EndKEY
\Address
Center for Geometry
{}~~~and Theoretical Physics
Box 90318
Duke University
Durham, NC 27708-0318
USA
\Address
\Address
\Address
\EndAddress

\section{Introduction}   \label{s:int}

Superstring theory does not currently have a complete definition. What
we have instead are a set of incomplete definitions each of which fill
in some of the unknown aspects of the other partial
definitions. Naturally two questions immediately arise given this
state of affairs:
\begin{enumerate}
\item Is each partial definition consistent with the others?
\item How completely do the partial definitions combine to define
string theory?
\end{enumerate}

Neither of these questions has yet to be answered and indeed both questions
appear to be quite deep. The first of these concerns the
subject of ``string duality''. Let us list the set of known
manifestations of string theory each of which leads to a partial
definition:
\begin{enumerate}
\item Type I open superstring
\item Type IIA superstring
\item Type IIB superstring
\item $E_8\times E_8$ heterotic string
\item $\spnh$ heterotic string
\item Eleven-dimensional supergravity (or ``M-theory'')
\end{enumerate}
The first five of these ``theories'' describe a string, which is
closed in all cases except the first, propagating in flat
ten-dimensional Minkowski space $\R^{9,1}$.
The last theory is more like that of a membrane propagating in
eleven-dimensional Minkowski space $\R^{10,1}$. (Note that many people
like to think of string theory as a manifestation of M-theory rather
than the other way around.)

Instead of using a completely flat Minkowski space, one may try to
``compactify'' these string theories by replacing the Minkowski space
by $X\times M$, where $X$ is some compact $(10-d)$-dimensional manifold (or
$(11-d)$-dimensional in the case of M-theory) and
$M\cong\R^{d-1,1}$. So long as all length scales of $X$ are large with
respect to any natural length scale intrinsic to the string theory, we
can see that $X\times M$ may approximate the original flat Minkowski
space. This is called the ``large radius limit'' of $X$.
One of the most fascinating aspects of string theory is that
frequently we may also make sense of compactifications when $X$ is
small, or contains a small subspace in some sense. An extreme case of
this is when $X$ is singular. In particular, $X$ need not be a
manifold in general.

The key ingredient to be able to analyze string theories on general
spaces, $X$, is {\em supersymmetry}. For our purposes we may simply
regard a supersymmetry as a spinor representation of the orthogonal
group of the Minkowski space in which the string theory lives. In
general a theory may have more than one supersymmetry in which case
the letter ``$N$'' is commonly used to denote this number. In the
above theories the type I and heterotic strings together with M-theory
each have $N=1$ while the type II strings correspond to $N=2$.

Upon compactification, the value of $N$ will change depending upon the
global holonomy of the Levi-Cevita connection of the tangent bundle of
$X$. The new supersymmetries of $M$ are constructed from the
components of the old spinor representations of the original Minkowski
space which are invariant under this holonomy. We will give some
examples of this process shortly.

The general rule is that the more supersymmetry one has, the more
tightly constrained the string theory is and the easier it is to
analyze away from the large radius limit. Note that this rule really
depends upon the total number of components of all the supersymmetries
and so a large $d$ has the same effect as a large $N$ (since $M$ is
has $d$ dimensions and so its spinor representation would have a
large dimension).

As well as constraining the string theory so that it may be more
easily analyzed, supersymmetry can be regarded as a coarse
classification of compactifications. A knowledge of $d$ and
$N$ provides a great deal of information about the resulting system.
Almost every possibility for $d$ and $N$ is worth at least one long
lecture in itself. We will deal with the case $d=4$ and $N=2$ about
which probably the most is known at this present time.

The principle of duality can now be stated as follows. Given a
specific string and its compactification on $X$ can
one find another string theory  compactified on another space, $Y$,
such that the ``physics'' in the uncompactified space, $M$, is
isomorphic between the two compactifications? This is important if our
first question of this
introduction is to be answered. In particular it should always be true
for any pair of string theories in our list unless there is a good
reason for a ``failure'' of one of the strings in some sense. We will
see an example of this below.

A mathematical analysis of duality requires a precise definition of
the physics of a compactification. This is not yet known in
generality. What we do know is a set of objects which are determined
by the physics, such as moduli spaces, partition functions,
correlation functions, BPS soliton spectra etc., which may be compared
to find necessary conditions for duality.

The most basic object one may study to identify the physics of two
dual theories is their moduli spaces. Roughly speaking this should
correspond to the moduli spaces of $X$ and $Y$ although one always
requires ``extra data'' beyond this. It is the extra data which leads
to the mathematical richness of the subject. Clearly if two theories are to be
identified, one must be able to identify their moduli spaces point by
point. This will be the focus of this talk.

It is a pleasure to thank my collaborators R.~Donagi and D.~Morrison
for many useful
discussions which were key to the results of section \ref{s:HII}.


\section{String Data}  \label{s:data}

In order to be able to describe the moduli space of each string theory
we are required to give the necessary data which goes into
constructing each one. Unfortunately, we do not have anywhere near enough
space to describe the origin of what follows. We refer to
\cite{Do:str,me:N2lect,GSW:book} for more details.

The theories which yield $d=4$ and $N=2$ in which we will interested are
specified by the following
\begin{itemize}
\item The type IIA string is compactified on a \CY\ threefold $X$
(which has $\SU(3)$ holonomy). The following data specifies the theory.
\begin{enumerate}
  \item A Ricci-flat metric on $X$.
  \item A $B$-field $\in H^2(X,\R/\Z)$.
  \item A Ramond-Ramond (RR) field $\in H^{\textrm{odd}}(X,\R/\Z)$.
  \item A dilaton$+$axion, $\Phi\in\C$.
\end{enumerate}
\item The type IIB string is compactified on a \CY\ threefold $Y$
(which also has $\SU(3)$ holonomy). The following data specifies the theory.
\begin{enumerate}
  \item A Ricci-flat metric on $Y$.
  \item A $B$-field $\in H^2(Y,\R/\Z)$.
  \item A Ramond-Ramond (RR) field $\in H^{\textrm{even}}(Y,\R/\Z)$.
  \item A dilaton$+$axion, $\Phi\in\C$.
\end{enumerate}
\item The $E_8\times E_8$ heterotic string is compactified on a
product of a K3 surface, $Z$, and an elliptic curve, $E_H$.
This product has $\SU(2)$ holonomy.
The following data specifies the theory.
\begin{enumerate}
  \item A Ricci-flat metric on $Z\times E_H$.
  \item A $B$-field $\in H^2(Z\times E_H,\R/\Z)$.
  \item A vector bundle $V\to(Z\times E_H)$ with a connection
  satisfying the Yang--Mills equations and whose structure group
  $\subseteq E_8\times E_8$. The respective characteristic classes in
  $H^4$ for $V$ and the tangent bundle of $Z\times E_H$ are fixed to
  be equal.
  \item A dilaton$+$axion, $\Phi\in\C$.
\end{enumerate}
\end{itemize}

In each case we can only expect the data to provide a faithful
coordinate system in some limit.
This is a consequence of the the fact that we only really have a
partial definition of each string theory.
A sufficient condition for faithfulness is that
the target space is large --- i.e., all minimal cycles have a large
volume, {\em and\/} $|\Phi|\gg1$. Beyond this we may expect ``quantum
corrections''. In general the global structure of the moduli space can
be quite incompatible with this parameterization  --- it is only
reliable near some boundary.

On general holonomy arguments (see, for example, \cite{CFG:II,me:lK3})
one can argue that the moduli space factorizes locally
\begin{equation}
  \M \cong \M_H \times \M_V,
\end{equation}
where (at least at smooth points) $\M_H$ is a quaternionic K\"ahler
manifold and $\M_V$ is a special K\"ahler manifold. We refer the
reader to \cite{Frd:SK} for the definition of a special K\"ahler
manifold. These restricted holonomy types are expected to remain {\em exact\/}
after quantum corrections have been taken into account.

We may now organize the above parameters into how they span $\M_H$ and
$\M_V$. First we note that Yau's theorem \cite{Yau:} tells us that the
Ricci-flat metric on a \CY\ manifold is uniquely determined by a choice of
complex structure and by fixing the cohomology class of the K\"ahler
form, $J\in H^2(\bullet,\R)$. We may combine $J$ and $B$ to form
the ``complexified K\"ahler form'' $B+iJ\in H^2(\bullet,\C/\Z)$. We
then organize as follows
\begin{itemize}
  \item The Type IIA string: $\M_V$ is parametrized by the complexified
  K\"ahler form of $X$. $H^{\textrm{odd}}(X,\R/\Z)\cong H^3(X,\R/\Z)$
  is the {\em intermediate Jacobian\/} of $X$ and is thus an abelian
  variety. We then expect a factorization
  $\M_H\cong\C\times\M_H'$, where $\Phi$ is the coordinate along the
  $\C$ factor. Finally we have a fibration $\M_H'\to\M_{\textrm{cx}}(X)$
  with fibre given by the intermediate Jacobian, and
  $\M_{\textrm{cx}}(X)$ is the moduli space of complex structures on
  $X$.
  \item The Type IIB string: $\M_V$ is now parametrized by the complex
  structure of $Y$. $H^{\textrm{even}}(Y,\R/\Z)\cong
  H^0(Y,\R/\Z)\oplus H^2(Y,\R/\Z)\oplus H^4(Y,\R/\Z)\oplus
  H^6(Y,\R/\Z)$ may be viewed as an abelian
  variety. We again expect a factorization
  $\M_H\cong\C\times\M_H'$, where $\Phi$ is the coordinate along the
  $\C$ factor. Finally we have a fibration $\M_H'\to\M_{\textrm{Kf}}(Y)$
  with fibre given by the RR fields, and
  $\M_{\textrm{Kf}}(Y)$ is the moduli space of the complexified
  K\"ahler form of $Y$.
  \item The $E_8\times E_8$ heterotic string: Let us first assume that
  the bundle $V\to(Z\times E_H)$ factorizes as $(V_Z\to Z)\times
  (V_E\to E_H)$. Thus the structure group of $V_Z$ times the structure
  group of $V_E$ is a subgroup of $E_8\times E_8$. We now expect
  $\M_V$ to factorize as $\C\times\M_V'$, where $\Phi$ is the
  coordinate along the
  $\C$ factor (see \cite{FvP:Ka} for a more precise
  statement). $\M_V'$ is then the total moduli space of $V_E\to E_H$
  including deformations of the complex structure and the complexified
  K\"ahler form of $E_H$. $\M_H$ is the total moduli space of the
  fibration $V_Z\to Z$ including deformations of the Ricci-flat metric
  of $Z$.
\end{itemize}

Again we emphasize that the above statements are approximate and only
valid when the
target space is large and $|\Phi|\gg1$. They should be exact only at the
boundary of the moduli space corresponding to these limits.
It is important to see that factorization of the moduli space will restrict the
way that the quantum corrections may act. For example, in the type IIA
string the dilaton, $\Phi$, lives in $\M_H$. This means that $\M_V$
cannot be subject to corrections related to having a finite
$|\Phi|$. Equally, the K\"ahler form parameter governs the size of $X$
and so $\M_H$ will not be subject to corrections due to finite size.

It is this property that some parts of the moduli space can be free
from quantum corrections and that the interpretation of this part can
vary from string theory to string theory which lies at the heart of
the power of string duality. If two theories are simultaneously exact
at some point in the moduli space then we may address the first
question in our introduction. If at every point in the moduli space
some theory (perhaps as yet unknown) is in some sense exact then we
may address the second question.


\section{Mirror Duality} \label{s:mir}

Mirror symmetry as first suggested in \cite{GP:orb,CLS:mir} was a
duality between ``conformal field theories''. We may make a different
version of mirror symmetry, a little more in the spirit of ``full'' string
theories, by proposing the following \cite{AM:Ud}:
\begin{defn}
The pair $(X,Y)$ of \CY\ threefolds is said to be a {\em mirror\/}
pair if and
only if the type IIA string compactified on $X$ is physically
equivalent to the type IIB string compactified on $Y$.
\end{defn}
Of course, this definition is mathematically somewhat unsatisfying as
it depends on physics. However, it encompasses previous definitions of
mirror symmetry. We also assume the following
\begin{prop}
If $(X,Y)$ is a mirror pair then so is $(Y,X)$.
\label{p:mr}
\end{prop}
While this proposal is obvious from the old definitions it is not
completely clear that we may establish it rigorously using the above
definition.

Applying this to the moduli space description in the previous section
we immediately see that, ignoring quantum corrections,
$\M_{\textrm{Kf}}(X)$ should be identified with $\M_{\textrm{cx}}(Y)$
and equally $\M_{\textrm{Kf}}(Y)$ should be identified with
$\M_{\textrm{cx}}(X)$. We know that
$\M_V$ is unaffected by $\Phi$ corrections and we expect
$\M_{\textrm{cx}}(Y)$ to be {\em exact\/} since it is also unaffected by size
corrections.

We expect that $\M_{\textrm{Kf}}(X)$ be affected by size
corrections. Similarly, given proposition \ref{p:mr}, $\M_{\textrm{cx}}(X)$ is
exact and $\M_{\textrm{Kf}}(Y)$ will suffer from size corrections. We
will use the notation $\mQ$ to refer to a fully corrected moduli
space. Thus $\mQ_{\textrm{Kf}}(X)\cong\mQ_{\textrm{cx}}(Y)\cong
\M_{\textrm{cx}}(Y)$ but $\mQ_{\textrm{Kf}}(X)\ncong\M_{\textrm{Kf}}(X)$.

The corrections to $\M_{\textrm{Kf}}(X)$ take the form of
``world-sheet'' instantons and were studied in detail in celebrated
work of Candelas et al \cite{CDGP:}. In particular, the assertion that
$\mQ_{\textrm{Kf}}(X)\cong\mQ_{\textrm{cx}}(Y)$ allows one to count
the numbers of rational curves on $X$. Subsequently
a great deal of work (see for example \cite{Kon:mir,RT:GW,Giv:mir,LLY:GW})
has been done which has made this curve counting much more rigorous.

As well as $\M_{\textrm{Kf}}$ and $\M_{\textrm{cx}}$, it is instructive
to look at the abelian fibres of $\M_H$ in the context of mirror
symmetry. The effect of equating $\M_{\textrm{cx}}(X)$ with
$\M_{\textrm{Kf}}(Y)$ is to equate
\begin{equation}
  H^3(X,\Z) \sim H^0(Y,\Z)\oplus H^2(Y,\Z)\oplus H^4(Y,\Z)\oplus
  H^6(Y,\Z),  \label{eq:decomp}
\end{equation}
but that we expect this correspondence to make sense only if $Y$ is very
large. Note that by going around closed loops in $\M_{\textrm{cx}}(X)$
we expect to have an action on $H^3(X,\Z)$ induced by monodromy. If we
were to take (\ref{eq:decomp}) to be literally true then we have to
say the same thing about the action of closed loops in
$\M_{\textrm{Kf}}(X)$ acting on the even integral cycles in $Y$. That
is to say, we would be claiming that if one begins with, say, a point
representing an element of $H^0(Y,\Z)$ we could smoothly shrink $Y$ down to
some small size and then smoothly let it re\"expand in some
inequivalent way such that our point had magically transformed itself
into, say, a 2-cycle! Clearly this does not happen in classical
geometry.

The suggestion therefore \cite{Cand:mir,AL:qag} is that quantum
corrections should be applied to the notion of integral cycles so
that, in the context of stringy geometry, 0-cycles {\em can\/} turn
into 2-cycles when $Y$ is small. Thus the notion of dimensionality
must be uncertain for small cycles.

Of central importance to the study of mirror pairs is being in a
region of moduli space where the quantum corrections are small. That
is we require $Y$ to be large. This amounts to a specification of the
K\"ahler form on $Y$ and must therefore specify some condition on the
complex structure of $X$. This was analyzed by Morrison:
\begin{prop}
  If $Y$ is at its large radius limit then $X$ is at a degeneration
  of complex structure corresponding to maximal unipotent monodromy.
\end{prop}
We refer the reader to \cite{Mor:gid} for an exact statement of this.
The idea is that $X$ degenerates such that a variation of mixed Hodge
Structures around this point leads to monodromy compatible with
(\ref{eq:decomp}).

The point we wish to emphasize here is that when $X$ is very large
then the complex structure of $Y$ is restricted to be very near a
particular point in $\M_{\textrm{cx}}(X)$. We only really expect
mirror symmetry to be ``classically'' true at this degeneration. Close to this
degeneration we may measure quantum perturbations leading to such
effects as counting rational curves. A long way from this degeneration
mirror symmetry is much more obscure from the point of view of classical
geometry.

It is possible to have a \CY\ threefold, $X$, whose moduli space
$\M_{\textrm{cx}}(X)$ contains no points of maximal unipotency. In
this case, its mirror, $Y$, can have no large radius limit. Since
clearly any classical \CY\ threefold may be taken to be any size, $Y$
cannot have an interpretation as a \CY\ threefold. This is the sense
in which duality can sometimes break down.


\section{Heterotic/Type IIA Duality}    \label{s:HII}

Having discussed mirror duality between the type IIA and the type IIB
string we will now try to repeat the above analysis for the duality
between the type IIA and the $E_8\times E_8$ heterotic
string. This duality was first suggested in \cite{KV:N=2,FHSV:N=2}
following the key work of \cite{HT:unity,W:dyn}.

In this case $\M_V$ is currently fairly well-understood (see,
for example \cite{me:lK3} and references therein). Here we will
discuss $\M_H$ which provides a much richer structure.

First let us discuss the quantum corrections.  On the heterotic side,
$\M_H$ contains the deformations of $Z$ as well as the vector bundle
over it. Note that in the the case of K3 surfaces we may not factorize
the moduli space of Ricci-flat metrics into a moduli space of complex
structures and the K\"ahler cone. This follows from the fact that
given a fixed Ricci-flat metric, we have an $S^2$ of complex
structures. The size of the K3 surface is a parameter
of $\M_H$ and so we expect $\M_H$ to suffer from quantum corrections
due to size effects for the heterotic string.

We also know that on the type IIA side, the dilaton is contained in
$\M_H$. Thus we expect $\M_H$ to suffer from corrections due to $\Phi$
for the type IIA string. We managed to evade worrying about such
effects in our discussion of mirror symmetry but here we are not so
lucky.

Let us now attempt to find the place in the moduli space where we may
ignore the quantum effects both due to $\Phi$ and due to size. To do
this we require the following:
\begin{prop}
  If a type IIA string compactified on a \CY\ threefold $X$ is dual to
  a heterotic string compactified on a factorized bundle over a
  product of a K3 surface, $Z$, and an elliptic curve $E_H$, then $X$
  must be in the form of an elliptic fibration $\pi_F:X\to\Sigma$ with
  a section and a K3 fibration $\pi_A:X\to B$. Here $\Sigma$ is a
  birationally ruled surface and $B\cong\P^1$.
\end{prop}
Note that these fibrations may contain degenerate fibres.
We refer to \cite{me:lK3} for details.

Let us now assume that $Z$ is in the form of an elliptic fibration
over $B$ with a section. Given this, we claim the following:
\begin{prop}
  The limit of large $Z$ automatically ensures that $\Phi\to\infty$
  for the type IIA string. In this limit, $X$ also undergoes a degeneration
  to $X_1\cup_{Z_*}X_2$, where $X_1$ and $X_2$ are each elliptic
  fibrations over a birationally ruled surface and are each fibrations
  over $B\cong\P^1$ with generic fibre given by a rational elliptic
  surface (RES). $Z_*= X_1\cap X_2$ is isomorphic to $Z$ as a complex variety.
\label{p:deg}
\end{prop}
We refer to \cite{FMW:F,AM:po} for a proof.

Recall that a RES is a complex surface given by
$\P^2$ blown up at nine points given by the intersection of two cubic
curves. In a sense, for elliptic fibrations a RES is ``half of a K3
surface''. This degeneration is viewed as each K3 fibre of the
fibration $\pi_A:X\to B$ breaking up into two RES's.

This degeneration therefore provides the analogue of the ``maximally
unipotent'' degeneration in the case of mirror symmetry. There are
important differences however. Note that while the
maximally unipotent degeneration of mirror symmetry essentially
corresponds to a point in the moduli space of complex structures, the
degeneration given by proposition \ref{p:deg} is not rigid --- it
corresponds a family of dual theories. In the case of mirror symmetry,
by taking $Y$ to be large we needed to fix a point in
$\M_{\textrm{Kf}}(Y)$ and thus $\M_{\textrm{cx}}(X)$. Here we need to
take the K3 surface $Z$ to its large radius limit but this does {\em
not\/} fix a point in $\M_H$. We may still vary the complex structure
of $Z$ (subject only to the constraint that it be an elliptic
fibration with a section) and we may still vary the bundle $V_Z$.

We should therefore be able to see the moduli space of complex
structures on the elliptic K3 surface, $Z$, and the moduli space of
the vector bundle $V_Z$ {\em exactly\/} from this degeneration of
$X$. The correspondence $Z\cong Z_*=X_1\cap X_2$ tells us how the
moduli space of $Z$ can be seen from the moduli space of the
degenerated $X$. The moduli space of the vector bundle is a little
more interesting.

$V_Z$ may be split into a sum of two bundles $V_{Z,1}$ and $V_{Z,2}$
each of which has a structure group $\subseteq E_8$. We will identify
$V_{Z,1}$ from a curve $C_1\subset Z_*$ and $V_{Z,2}$ from $C_2\subset Z_*$.
$C_1$ and $C_2$ will form the {\em spectral curves\/} of their
respective bundles in the sense of \cite{Don:spec}.

Let us consider a single RES fibre $Q_b$ of the fibration $X_1\to B$,
where $b\in B$. $Q_b$ is itself an elliptic fibration $\pi_Q:Q_b\to\P^1$.
The section of the elliptic fibration $\pi_F:X\to\Sigma$ determines a
distinguished section $\sigma_0\subset Q_b$.
Blowing this down gives a Del Pezzo surface with 240
lines $\sigma_1,\ldots,\sigma_{240}$.

We then have
\begin{prop}
  The fibre of the branched cover $C_1\to B$ is given by the set
  of points $\{\sigma_i\cap Z_*; i=1,\ldots,240\}$,
\end{prop}
with an analogous construction for $C_2$. We refer to
\cite{me:hyp,CD:F4} for details.

We also have the data from the abelian fibre of $\M_H$ corresponding
to the RR fields. In the case of heterotic/type IIA duality
we have \cite{me:hyp}
\begin{prop}
$$
  \Lambda_0 \cong H^1(C_1,\Z) \oplus H^1(C_2,\Z) \oplus H^2_T(Z,\Z),
$$
where $\Lambda_0$ is the sublattice of $H^3(X,\Z)$ invariant under
monodromy around the degeneration of proposition \ref{p:deg} and
$H^2_T(Z,\Z)$ is the lattice of transcendental 2-cocycles in $Z$.
\label{p:HetM}
\end{prop}
Thus the RR-fields of the type IIA string map to the Jacobians of $C_1$
and $C_2$, required to specify the bundle data, and to the $B$-field
on $Z$.

Proposition \ref{p:HetM} should embody much of the spirit of the
duality between the type IIA string and the $E_8\times E_8$ heterotic
string in a similar way that equation (\ref{eq:decomp}) embodies
mirror symmetry. In particular $\Lambda_0$ is not invariant under
monodromy around any loop in the moduli space and so the notion of
what constitutes the $E_8$-bundles and what constitutes the K3 surface
$Z$ should be blurred in general --- just as the notion of 0-cycles
and 2-cycles is blurred in mirror symmetry.

The analysis of the moduli space $\M_H$ is very much in its
infancy. In this talk we have not even mentioned how to compute
quantum corrections --- the above discussion was purely for the exact
classical limit. There appear to be many adventures yet to be encountered
in bringing the understanding of heterotic/type IIA duality to the
same level as that of mirror symmetry.


\Addresses
\end{document}